\documentclass{IEEEtran4PSCC}
\usepackage{hyperref}
\usepackage{nameref}
\usepackage{varioref}
\usepackage{graphicx}
\usepackage{color}
\usepackage{bigints}
\usepackage{amssymb}
\usepackage{amsmath}
\usepackage{dsfont}
\usepackage{mathtools}
\usepackage{comment}
\usepackage{cleveref}
\usepackage{multirow}
\usepackage{amsfonts}
\usepackage{latexsym}
\usepackage{algorithm,algpseudocode}
\usepackage{ctable}

\let\oldReturn\Return
\renewcommand{\Return}{\State\oldReturn}
\usepackage{enumerate}
\usepackage{wasysym}
\usepackage{lipsum}
\usepackage{amssymb}
\usepackage{amsthm,bm}
\usepackage{listings}
\usepackage{color}
\usepackage{graphicx}
\usepackage[font=small,labelfont=bf]{caption} 
\usepackage[T1]{fontenc}
\usepackage{beramono}
\usepackage{listings}
\usepackage{amsmath}

\usepackage{xurl}
\definecolor{mustard}{HTML}{707238}
\definecolor{PaleGreen}{HTML}{98fb98}
\definecolor{LightRed}{HTML}{fe7a7a}
\definecolor{LightYellow}{HTML}{effa86}
\definecolor{YellowOrange}{HTML}{ffae42}
\definecolor{OliveGreen}{HTML}{408000}
\definecolor{dkgreen}{rgb}{0,0.6,0}
\definecolor{gray}{rgb}{0.5,0.5,0.5}
\definecolor{mauve}{rgb}{0.58,0,0.82}
\definecolor{darkred}{rgb}{139,0,0}
\title{Real-Time Stochastic Assessment of Dynamic N-1 Grid Contingencies}
\author{
\IEEEauthorblockN{Ayrton Almada \IEEEauthorrefmark{1}\IEEEauthorrefmark{2}, Laurent Pagnier \IEEEauthorrefmark{1}, Igal Goldshtein \IEEEauthorrefmark{3}, Saif R. Kazi \IEEEauthorrefmark{2},
and Michael (Misha) Chertkov \IEEEauthorrefmark{1}}
\IEEEauthorblockA{\IEEEauthorrefmark{1} Program in Applied Mathematics and Department of Mathematics,\\
University of Arizona, Tucson, AZ 85721, USA}
\IEEEauthorblockA{ \IEEEauthorrefmark{2} Applied Mathematics and Plasma Physics,\\ Los Alamos National Laboratory, Los Alamos, NM 87519, USA}
\IEEEauthorblockA{\IEEEauthorrefmark{3} NOGA, Power System Operator of Israel, Haifa, Israel}
}
\makeatletter
\let\old@ps@headings\ps@headings
\let\old@ps@IEEEtitlepagestyle\ps@IEEEtitlepagestyle
\def\psccfooter#1{
    \def\ps@headings{
        \old@ps@headings
        \def\@oddfoot{\strut\hfill#1\hfill\strut}
        \def\@evenfoot{\strut\hfill#1\hfill\strut}
    }
    \def\ps@IEEEtitlepagestyle{
        \old@ps@IEEEtitlepagestyle
        \def\@oddfoot{\strut\hfill#1\hfill\strut}
        \def\@evenfoot{\strut\hfill#1\hfill\strut}
    }
    \ps@headings
}
\makeatother
\psccfooter{}
\begin{document}
\maketitle
\begin{abstract}
Power system operators need tools for rapid, real-time counterfactual assessments of grid security under fast-changing conditions. Traditional N-1 contingency analysis lacks dynamic evaluation, especially of frequency swings from common faults. This paper introduces a real-time dashboard framework to screen dynamic contingencies. It assumes: (a) the grid starts in a balanced state; (b) faults can occur randomly on any transmission line, temporarily de-energizing and then reconnecting it within about one second; and (c) contingencies are flagged if post-fault transients cause line flows to exceed safety thresholds. The key contributions are: (1) Overload Indicator – a system-wide metric quantifying integrated N-1 dynamic risk from a given state; (2) Scalable Fault Evaluation Algorithm – a linear-scaling method to assess dynamic fault impacts without brute-force simulations; and (3) Risk Estimation – a Cross Entropy Adaptive Importance Sampling method estimating the likelihood of low probability by high risk events, e.g. associated with potential transformer over-current. We demonstrate the framework on the Israeli transmission power grid (IG).
\end{abstract}
\thanksto{\noindent The work at UA was supported by (a) NSF DMS-2229012: "Collaborative Research: AMPS: Rare Events in Power Systems: Novel Mathematics, Statistics and Algorithms"; (b) subcontract from NOGA on "Uncertainty-Aware Toolbox for Simulation, Optimization, Control and Planning of Inter Connected Gas and Power Systems". AA was supported by LANL as a graduate research assistant. We are thankful to Prof. J. Blanchet (Stanford) and his team for multiple discussions and consultations on the application of the adaptive importance sampling (cross-entropy) method; and Dr. Bent (LANL) for multiple comments and suggestions.}
\section{Introduction}
\label{sec:intro}
\subsection{Motivation}
Power-system stability is the grid’s ability to regain an acceptable operating state after a disturbance without loss of system integrity. It ensures that perturbations -- from routine load/generation fluctuations to severe faults -- do not trigger unnecessary tripping or cascading outages. Stability margins depend on network topology, the pre-disturbance operating point, disturbance characteristics, and protection/control settings; sustained synchronous operation of electromechanical modes is central to this resilience \cite{machowski2020power}. Transient stability focuses on maintaining synchronism in the milliseconds-to-tens-of-seconds window immediately after a disturbance and is governed by the initial operating state, the disturbance severity, and the calibration/coordination of automatic protection. With growing system heterogeneity and uncertainty (e.g., variable renewables, demand response, advanced power electronics), {\bf real-time dynamic-security assessment} is increasingly necessary, complementing -- rather than replacing -- traditional offline studies \cite{xue2006power}. Recent continent-scale and U.S. incidents illustrate how fast voltage/frequency and protection interactions can tip the grid into low-probability, high-impact outcomes (Iberia, 28 April 2025; Texas, February 2021) \cite{entsoe2025iberian,fercnerc2021uri}.

In this manuscript we examine the possibility of creating such a new on-line tool -- a {\bf dashboard to inform the system operator about counterfactual (possible) contingencies}. Specifically, we are concerned about resolution of faults, e.g. single- and three- phase faults at any (of many) lines in the system, and thus estimate probability of a secondary effect -- rare but potentially extremely damaging -- caused by such a fault -- with one of a critical elements of the system, e.g. a transformer, becoming overload during the primary post-fault seconds-long transient.

\subsection{Main Contributions}

\begin{figure}[thb]
    \centering
    \includegraphics[width=1.0\columnwidth]{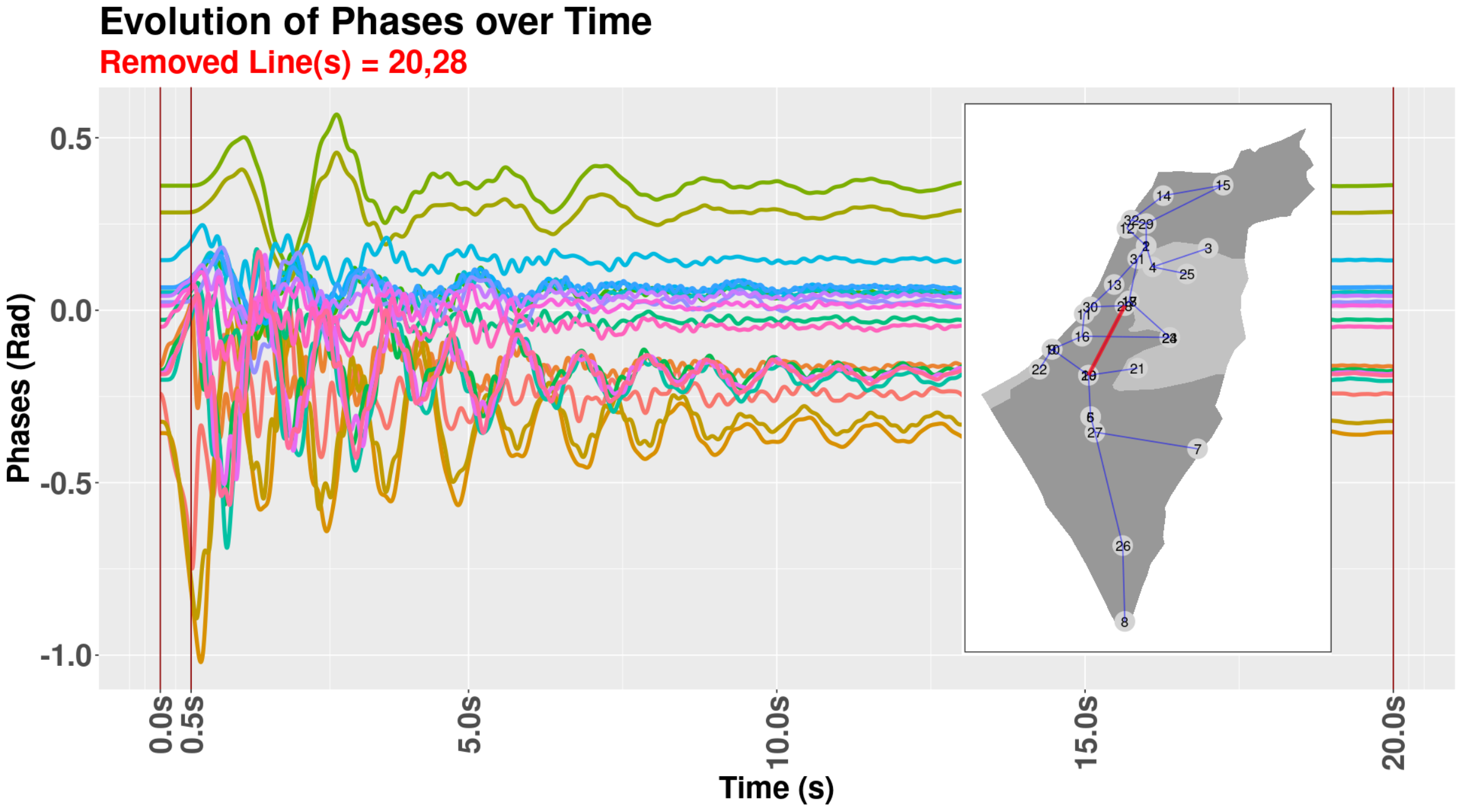}
	\caption{Single-Phase Fault Example (IG) – Unstable Case. If the faulty line isn’t re-energized quickly, loss of synchronization occurs (system splits). Fault, Clearance, Finish occurs at time: 0.0, 0.5 and 20s respectively.
    } 
	\label{Point3}
\end{figure}

This project introduces 
\begin{itemize}
    \item \textbf{An efficient methodology for $N-1$ dynamic screening} -- analyze and simulate the dynamics of an ensemble of transients associated with counterfactual faults of any line for a given state. Fig.~(\ref{Point3}) shows an example of a undesirable transient behavior that results from triggering a fault in the system that lasts $0.5s$.
    \item \textbf{A probabilistic framework} to estimate how system overloads evolve over time based on angular phase behavior, supporting assessment of grid reliability under stochastic disturbances.
\end{itemize}

We obtain a \textbf{three–orders-of-magnitude} (i.e., $\sim 10^3$) efficiency gain over direct Monte Carlo of counterfactual contingencies via three cumulative improvements:
\begin{enumerate}
  \item \textbf{Analytic reduction:} replace direct time–domain simulations with linear-algebraic solutions of the governing models;
  \item \textbf{Iterative reuse:} approximate the analytic solution with fast iterations that \emph{reuse} factorizations/pre-conditioners across many contingencies, thereby avoiding a full large matrix inversion per case;
  \item \textbf{Targeted sampling:} employ cross-entropy–based importance sampling to concentrate trials on high-impact contingency candidates.
\end{enumerate}

\subsection{Outline}

The rest of this paper is divided into the following sections \textbf{\ref{sec:model}}-\textbf{\ref{sec:conclusion}}:

\begin{itemize}
    \item \textbf{System and Fault Modeling} [Section \textbf{\ref{sec:model}}]:\\ Develop a scalable, linearized second-order model of power grid dynamics based on swing equations. Model both single- and three-phase faults to simulate disturbances over three intervals: pre-fault, faulted, and post-clearing.
    \item \textbf{Overload Indicator} [Section \textbf{\ref{sec:indicator}}]:\\ Introduce a metric that captures fault severity using thermal limits, fault duration, and phase oscillations—providing a continuous assessment of system response.
    \item \textbf{Scalable Dynamic Evaluation} [Section \textbf{\ref{sec:algorithm}}]:\\ Apply analytics to re-solve the system's dynamics efficiently.
    \item \textbf{Probabilistic Risk Estimation} [Section \textbf{\ref{sec:mcmc}}]:\\ Estimate rare-event probabilities using Monte Carlo and advanced importance sampling method --  efficient Cross-Entropy (validated on brute force MCMC) --  through our custom simulator, \textbf{N1Plus}.
    \item \textbf{Case Study: Israeli Grid}  [Section \textbf{\ref{sec:casestudy}}]:\\ Demonstrate the framework’s accuracy and efficiency on the IG network (32 buses, 36 lines, 11 generators) under various fault scenarios.
    \item \textbf{Conclusion and Outlook} [Section \textbf{\ref{sec:conclusion}}]:\\ Extend the framework to more complex grids, e.g. utilizing AI tools to scale it up. 
\end{itemize}

We illustrate and summarize the constitution of this work on Fig.~(\ref{fig:pipeline}):
\begin{figure}[thb]
  \centering
  \includegraphics[width=1.0\columnwidth]{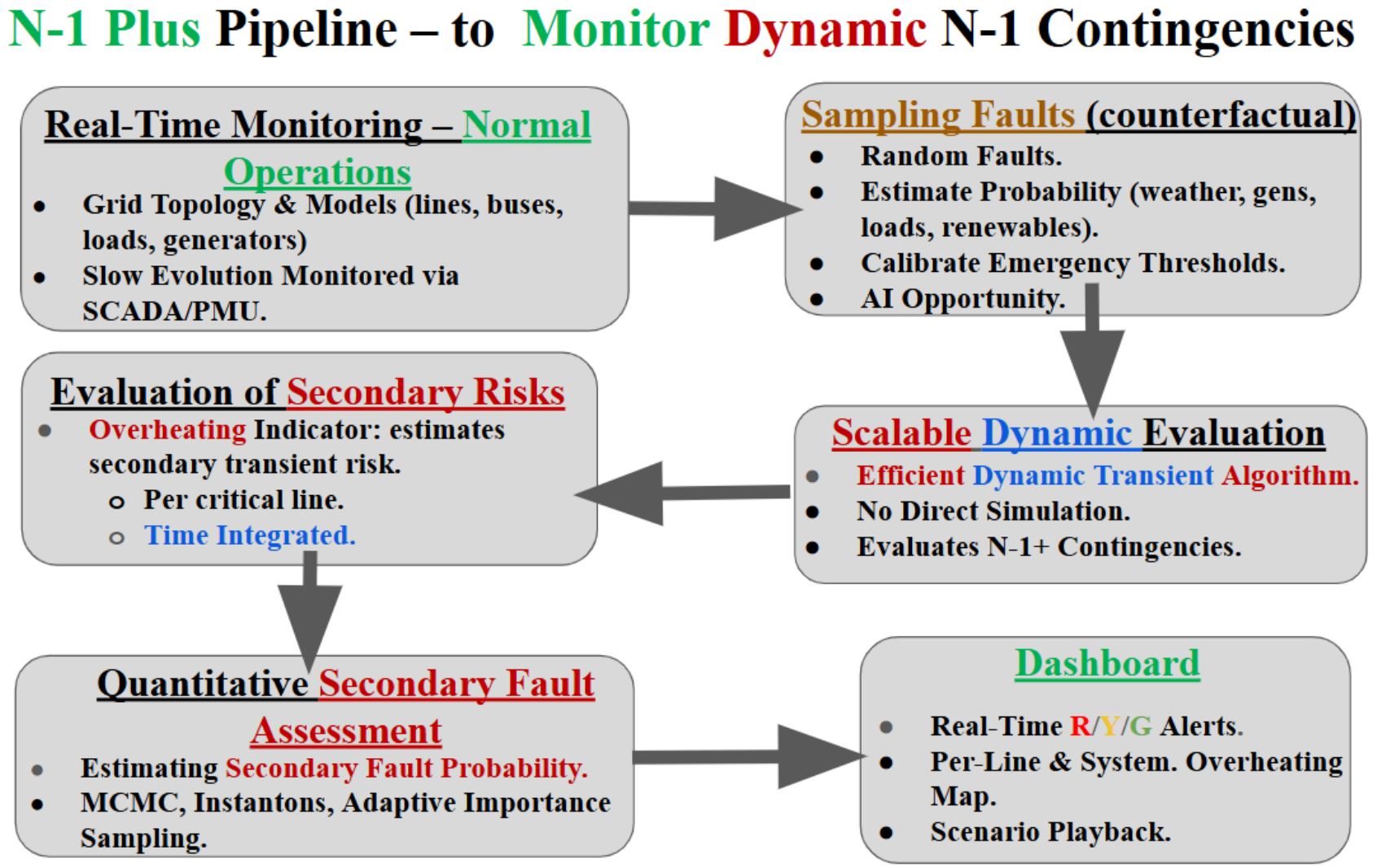}
  \caption{Conceptual Dashboard for Dynamic Contingency Screening.}
  \label{fig:pipeline}
\end{figure}

\section{Related Work and Background}
\label{sec:related}

\subsection{$N-1$ Contingency Analysis}

The $N-1$ contingency analysis\footnote{Our methodology also generalizes to broader classes of contingencies, including those defined from historical data rather than strict $N-1$ rules.} is a fundamental reliability tool in power system planning and operation. It ensures that the grid remains stable and within operational limits after the failure of any single critical component, such as a line, generator, or transformer. This principle underlies the $N-1$ security standard adopted by organizations such as NERC and ENTSO-E \cite{chatterjee2010n}. Traditionally, $N-1$ analysis simulates individual outages from a base-case power flow to identify violations (e.g., thermal or voltage limits) and flag critical contingencies for corrective actions such as re-dispatch or network reconfiguration \cite{milano2010power}\footnote{Voltage regulation and dispatch corrections are not addressed in this manuscript.}. While effective, these static methods are computationally demanding and often overlook dynamic effects, especially in large, renewable-rich systems with low inertia and high variability. To address these limitations, we introduce a framework for real-time screening of \textit{dynamic} $N-1$ contingencies, focusing on short-term transient behavior (sub-second to tens of seconds; see Table~\ref{tab:litreview}) following typical disturbances such as single-phase short circuits. Our method estimates the probability of secondary risks -- seconds-long but severe overloads on critical lines -- that may trigger further instabilities in an already stressed grid. The ultimate goal is to provide operators with a clear, actionable dashboard that enhances situational awareness of rare but potentially severe post-fault events, such as line tripping induced by transient overcurrents.

\begin{table*}[!t]
\caption{N-1 Analysis vs Dynamic Evaluation.}
  \label{tab:litreview}
\small
\centering
\begin{tabular}{lll}
\hline
\multicolumn{1}{c}{\textbf{Aspect}}
 & \multicolumn{1}{c}{\textbf{N-1 Contingency}} & \multicolumn{1}{c}{\textbf{Swing Equations}} \vphantom{$\int_f^f$} \\
\hline
\vphantom{$\int_f^f$}Objective & Secure post-fault & Assess transient stability \vspace{1pt}\\
Time Scale [s] & $10^0 - 10^2$ & $10^{-3} - 10^{0}$ \vspace{3pt}\\
Application & Operations and planning& Dynamic stability  studies\vspace{1pt}\\
Model & AC/DC power flow &  Generator and load  responses as ODEs\vspace{1pt}\\
\hline
\end{tabular}
\end{table*}

\subsection{Rare Event/Instanton Sampling}

Rare–event (instanton) methods provide an efficient framework for identifying low-probability, high-impact events in complex systems \cite{juneja2006rare,asmussen2007stochastic,touchette2009large,rubino2009rare,freidlin2012random,hartmann2012efficient,grafke2015instanton}, including power grids \cite{dobson_complex_2007,chertkov2010predicting,chertkov_exact_2011,pfitzner_statistical_2011,owen2019importance}. Unlike standard Monte Carlo, these approaches leverage large-deviation theory to characterize the most probable rare trajectories -- \emph{instantons} -- that drive extreme outcomes such as blackouts and cascading failures \cite{dobson_complex_2007,pfitzner_statistical_2011}. The instanton minimizes a stochastic action (Freidlin–Wentzell) and thereby reveals both the physical mechanism and the exponential scaling of the associated probabilities \cite{juneja2006rare,freidlin2012random,grafke2015instanton}. Adaptive importance sampling (AIS), and in particular the cross-entropy (CE) method \cite{rubinstein_optimization_1997,de2005tutorial}, can be interpreted as using instanton-informed proposal families to bias sampling toward the most relevant regions of state space; iterative CE/AIS updates refine this bias via stochastic optimization, yielding substantial variance reduction and, in the low-probability regime, asymptotically optimal (near zero-variance) estimators \cite{rubinstein_optimization_1997,de2005tutorial,asmussen2007stochastic,hartmann2012efficient,owen_monte-carlo_2013}.

\section{Dynamics, Faults and Overloads}
\label{sec:model}

This section introduces the core modeling framework used throughout the paper: a transmission-level power grid dynamics model based on swing equations. We begin by formulating the system's behavior under small-to-moderate disturbances, then describe how line faults are incorporated into the model, and finally define a safety domain that encodes grid element overload constraints. This framework enables both theoretical analysis and efficient numerical simulation of dynamic contingencies.

\subsection{Modeling with Swing Equations}

A common power transmission grid model consists of generators (controlling resources), loads (including negative loads for renewable generators which are not in control), and the transmission lines interconnecting them. A generator bus includes both a synchronous generator and an associated load, while a pure load bus does not contain generation. Aggregated loads at the transmission level may exhibit their own dynamics -- typically on the order of minutes or slower. In contrast, generator dynamics, which are driven by mechanical inertia and respond to exogenous perturbations (e.g., changes in load, renewable injections, or fault-induced disturbances), evolve over shorter timescales ranging from sub-seconds to tens of seconds. This is the temporal regime of interest in our analysis. 

Each bus \( i \in \mathcal{V} \) is characterized by a complex voltage potential \( v_i = V_i e^{j\theta_i} \), where \( V_i \) is the voltage magnitude, \( \theta_i \) is the voltage phase angle, and \( j^2 = -1 \). In practice, transmission lines are three-phase systems; however, under the standard assumption of phase balance, we model the grid using a single equivalent phase per line \cite{bergen2009power}. 

Mathematically, the transmission grid is represented as an undirected graph \( \mathcal{G} = (\mathcal{V}, \mathcal{E}) \), where \( \mathcal{V} = \{1,2,\dots,|\mathcal{V}|\} \) denotes the set of buses and \( \mathcal{E} \subseteq \mathcal{V} \times \mathcal{V} \) the set of transmission lines. The set of generator buses is denoted \( G = \{1,2,\dots,|G|\} \), and load buses are labeled \( L = \{|G|+1, \dots, |\mathcal{V}|\} \). Throughout this work, we assume a lossless grid with constant voltage magnitudes \( V_k \) for all \( k \in \mathcal{V} \), appropriate for the sub-second to tens-of-seconds regime where voltage regulation and resistive losses are negligible.

The dynamical behavior of the system under small-to-moderate perturbations -- and these are regimes we consider -- is governed by the classical swing equations \cite{machowski2020power} (ignoring losses): 
\begin{align}
\label{eq:swing-nl}
\forall i \in \mathcal{V}:\quad m_i \ddot{\theta}_i \!+\! d_i \dot{\theta}_i 
\!+\!\!\!\sum_{\{i,j\}\in \mathcal{E}}\!\!\! V_iV_j B_{ij} \sin(\theta_i \!-\! \theta_j) \!=\! P_i,
\end{align}
where \( m_i \) denotes the inertia at bus \( i \), \( d_i \) the damping coefficient, \( P_i \) the net power injection (positive for generation, negative for consumption), and \( B_{ij} \) the susceptance of the transmission line connecting buses \( i \) and \( j \).

In typical operating conditions, disturbances trigger transients that dissipate within a few seconds. Once the system stabilizes, time derivatives vanish (\( \dot{\theta}_i = \ddot{\theta}_i = 0 \)), and Eq.~\eqref{eq:swing-nl} reduces to the static Power Flow (PF) equations.

Assuming the phase angle differences remain small during transients (i.e., \( |\theta_i - \theta_j| \ll 1 \)), we linearize the swing equations using the approximation \( \sin(\theta_i - \theta_j) \approx \theta_i - \theta_j \). This yields the following linearized swing model, which forms the foundation for the remainder of our analysis: 
\begin{align}
\label{eq:swing-lin}
m_i &\ddot{\theta}_i + d_i \dot{\theta}_i 
+ \sum_{\{i,j\} \in \mathcal{E}} \beta_{ij} (\theta_i - \theta_j) = p_i : i \in \mathcal{V}.
\end{align}
We introduce the vector form of Eq.~\eqref{eq:swing-lin} that will be used in Section \textbf{\ref{sec:algorithm}}:
\begin{align}\label{irtrp}
{\bm M} &\ddot{\theta} + {\bm D} \dot{\theta} 
+ {\bm L}\theta = p,
\end{align}
where ${\bm M},{\bm D}$ are the diagonal matrices for inertia and damping, ${\bm L}=\begin{pmatrix}L_{ij}\end{pmatrix}_{i\in\cal V}^{j\in\cal V}:$
$$
L_{ij}=\begin{cases}
\sum_{j\in\cal V}\beta_{ij}&:i=j\,,\\
-\beta_{ij}&:i\neq j\,.
\end{cases}
$$
is the \textbf{weighted graph Laplacian} of $\cal G$ and \( \beta_{ij} = V_iV_j B_{ij} \) represents the effective (scaled) \textbf{stiffness} of the line \( \{i,j\} \) \footnote{In general, the system of Eqs.~(\ref{irtrp}) should be viewed as a linearization of more general dynamic equations accounting for seconds-scale transients -- e.g. related to dynamics of generators/loads, inverters and transformers, possibly including effects of voltage -- with the inertia, ${\bm M}$-, damping, ${\bm D}$-, and graph-Laplacian, ${\bm L}$- matrices and the generator/load $p$- vector dependent on the current/base solution of the static power flow equations.}. Importantly, \( \beta_{ij} \) may be modified during a fault event -- either partially or entirely -- reflecting changes in the line's electrical properties. The modeling of such fault-induced changes is the focus of the next subsection.
\subsection{Simulating Faults}\label{sec:faults}
\begin{figure}[thb]
  \centering
  \includegraphics[width=1.0\columnwidth]{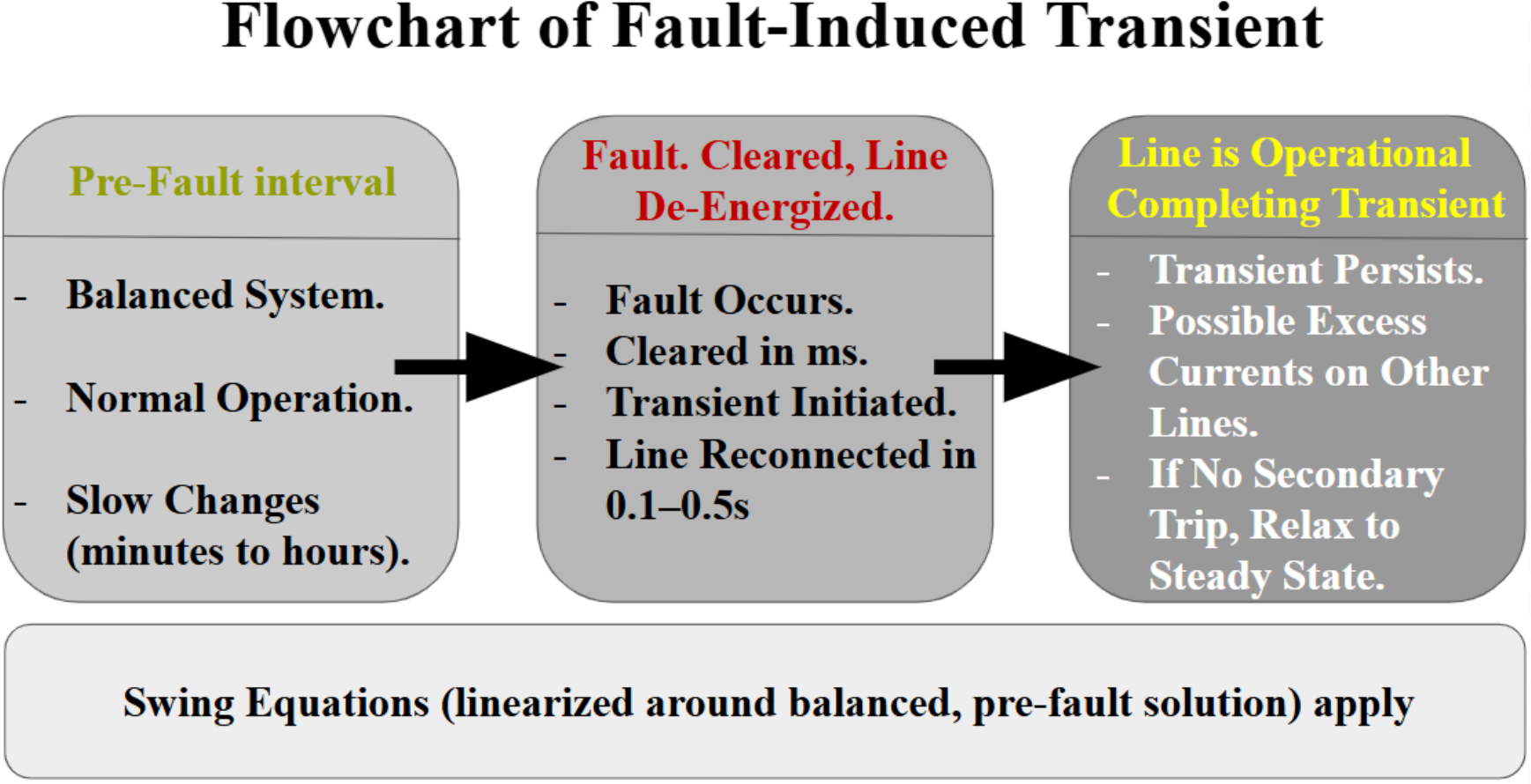}
  \caption{Fault timeline: system is balanced, pre-fault; fault is cleared, a faulty line is de-energized; line is back in service, post-fault transient. 
  }
  \label{fig:timeline}
\end{figure}

This study focuses exclusively on faults occurring on transmission lines. Faults at buses (e.g., transformers or generators) are excluded, although the proposed framework can be readily extended to include them. Two main fault types are considered: three-phase and single-phase faults\footnote{Intermediate fault severities could also be modeled, forming a continuous spectrum between single- and three-phase faults; such extensions are not explored here.}. Faults are assumed to occur randomly, with probabilities varying across lines according to seasonal conditions, weather, or system stress.
\textbf{Three-Phase Faults} represent the most severe disturbances in transmission networks. They may result from lightning, vegetation contact, or internal control actions enforcing symmetry during fast transients, and serve as a \textbf{benchmark for protection-device sizing and stability assessment} \cite{grainger1999power}. When a three-phase fault occurs on line $\{i,j\}$, the line is de-energized for a duration $\tau$ (typically tens to hundreds of milliseconds). In the model, this corresponds to removing the line by setting $\beta_{ij}=0$, i.e., \( \mathcal{E} \rightarrow \mathcal{E}_f = \mathcal{E} \setminus \{i,j\} \). After $\tau$, the line is restored, returning the network to its nominal topology. The duration $\tau$ may be deterministic—defined by relay settings—or stochastic, modeled as exponentially distributed with mean $\bar{\tau}$, ensuring oscillations remain below the nominal frequency limit (\textbf{50 Hz for the IG network}).
\textbf{Single-Phase-to-Ground Faults} are the most common, accounting for $70–80\%$ of all transmission disturbances \cite{horowitz2022power}. In many systems, regulations permit continued operation of the remaining two phases during such events. Under a balanced approximation, a single-phase fault is modeled by reducing the affected line’s susceptance from $\beta_{ij}$ to $\frac{2}{3}\beta_{ij}$ during the fault duration $\tau$, reflecting reduced transfer capability. Once cleared, $\beta_{ij}$ instantly returns to its nominal value\footnote{This simplification neglects unbalanced dynamics. Future work will incorporate more realistic single-phase fault models with unbalanced dynamic simulations.\label{footnote:disclaimer}}. As with three-phase faults, $\tau$ may be treated as deterministic or stochastic (e.g., exponentially distributed).
Fig.~(\ref{fig:timeline}) summarizes the three main stages in the system dynamics: \textcolor{mustard}{Pre-Fault}, \textcolor{red}{Fault}, and \textcolor{YellowOrange}{Operational} intervals.

\subsection{Safety Polytope}

We define the \emph{safety polytope} as the domain of phase configurations that ensure all line flows remain within their respective thermal (or angular) limits. Specifically, it characterizes the set of phase angles for which the power transmitted along every line does not exceed its safety threshold:

$$
\Pi_{\bm{\theta}} \!\doteq\! \left\{\! \bm{\theta} \!\in\! [-\pi,\pi]^{|\mathcal{V}|} \middle| \forall \{i,j\} \!\in\! \mathcal{E}:\left|\beta_{ij}\left(\theta_i \!-\! \theta_j\right)\right| \!\leq\! \bar{p}_{ij} \!\right\},
$$

where \( \bar{p}_{ij} \) denotes the maximum allowable phase difference across line \( \{i,j\} \), which indirectly encodes the corresponding thermal or stability constraint on power flow. The set \( \Pi_{\bm{\theta}} \) thus models the domain of safe operations for the system under the linearized swing approximation.

\section{Overload Indicator} 
\label{sec:indicator}

To quantify the \textbf{severity of line overloads} during dynamic transients, we define an \emph{overload indicator} based on the solution to the linearized swing equations, Eq.~\eqref{eq:swing-lin}. Let \(\theta(t) = (\theta_i(t)|i \in \mathcal{V})\) denote the vector of voltage phase angles over the observation interval \([0, T]\), and let \(\beta_{ij}(t)\) represent the time-dependent effective susceptance of line \(\{i,j\}\), as modeled throughout the fault (see Section~\ref{sec:faults}).

We define the \textbf{line-specific overload indicator} for a monitored line \(\{i, j\} \in \mathcal{E}\) as:
\begin{gather}
\label{eq:over-heat-ij}
S_{ij}(\theta_{0\to T})\! \doteq\!\! \int\limits_0^T\!\!\mathbb{I} \Big(\! \Big| \beta_{ij}(t)\big(\theta_i(t) \!-\! \theta_j(t)\big)\! \Big| \!>\! \bar{p}_{ij} \!\Big){\rm d}t ,
\end{gather}

where \(\mathbb{I}(\cdot)\) is the indicator function, returning 1 when the argument is true (i.e., the line is overloaded) and 0 otherwise. Here, \(\bar{p}_{ij}\) denotes the safety threshold for the phase difference (related to thermal limits) on line \(\{i,j\}\).

Importantly, the monitored line \(\{i,j\}\) is not necessarily the line where the fault originates. Since the location of the fault is unknown a priori, the indicator must account for dynamic overloads that arise elsewhere in the network as a consequence of any possible fault scenario.

To assess system-wide impact, we also define the \textbf{global overload indicator}:
\begin{gather}\label{glonOHI}
S(\theta_{0\to T}) \doteq \sum_{\{i,j\} \in \mathcal{E}_{\rm m}} S_{ij}(\theta_{0\to T}),
\end{gather}

where \(\mathcal{E}_{\rm m} \subseteq \mathcal{E}\) denotes the set of monitored lines -- typically selected based on criticality, though in smaller systems we may set \(\mathcal{E}_{\rm m} = \mathcal{E}\); and $\theta_{0\to T}\doteq (\theta_i(t)|i\in{\cal V}, t\in [0,T])$.

The total indicator \(S(\theta_{0\to T})\) accumulates the duration of overloads across monitored lines. Since the integrand in Eq.~\eqref{eq:over-heat-ij} is piecewise constant and integer-valued, the function grows in discrete steps -- each step corresponding to a time interval during which one or more lines exceed their limits. Thus, the indicator reflects both the \emph{extent} (number of overloaded lines) and the \emph{duration} of transient stress.

The quantities $S_{ij}(\theta_{0\to T})$ and $S(\theta_{0\to T})$ are the principal risk metrics evaluated across fault scenarios in this study. In a statistical setting, these indicators are computed across a distribution of fault events and system conditions (e.g., power injections), and their probability distributions are analyzed to assess grid vulnerability.

\begin{table*}[!t]
\caption{Comparison of Analytical vs Numerical Methods in number of bus $n$, time step $\Delta t$ and number of time steps $n_{\Delta t}\doteq T /\Delta t$.
}
  \label{Theovssim}
\small
\centering
\begin{tabular}{lll}
\hline
\multicolumn{1}{c}{\textbf{Feature}}
& \multicolumn{1}{c}{\textbf{Through Diagonalization}}  & \multicolumn{1}{c}{\textbf{Numerical Simulation }}\vphantom{$\int_f^f$}\\
\hline
\vphantom{$\int_f^f$}Computational cost per $X(t)$  eval.& $\mathcal{O}(n^3)$ twice, then $\mathcal{O}(n^2)$ & $\mathcal{O} \left( n_{\Delta t} \, n^2 \right)$ \\
Accuracy & Exact  & Depends on $\Delta t$\vspace{1pt}\\
\textbf{Evaluation at Multiple Time Points} & \textbf{Efficient after decomposition:
}  & \textbf{Requires full recomputation} \\
&$\mathcal{O}(n^2)$ per evaluation & for each time point\vspace{1pt}\\
\hline
\end{tabular}
\end{table*}

\section{Scalable Dynamic Evaluation} 
\label{sec:algorithm}

Linear swing Eqs.~(\ref{eq:swing-lin},\ref{irtrp}), describing  dynamics of phases over grids from time $0$ when the fault occurred and line (phase) was de-energized, through the period of time $\tau$ when the line as back in service, but swings of the phase dynamics continue, and all the way till time $T$ after the fault when the transient stabilized, can be restated in the following compact form for $t\in[0,T]$:
\begin{align}\label{eq:swing-lin-matrix}
    & \dot{\bm x}\!=\!\bm A(t)\,  \bm x + \bm P\,,\\ \nonumber &  \bm P\!\doteq\! \begin{pmatrix} \bm p \\ \bm 0\end{pmatrix},\ \bm x \!\doteq\! \begin{pmatrix} \dot{\bm \theta}(t) \\ \bm \theta(t)\end{pmatrix},\ \bm x_0\!\doteq\!\begin{pmatrix} \bm 0\\
    \bm L_0^{-1} \bm p
    \end{pmatrix},\\
    &\bm A(t) = \left[\!\!\begin{array}{cc}
    \bm M^{-1} \bm D & \bm M^{-1} \bm L(t)\\
    \mathbb{O}_{n\times n} & \mathds{1}_{n\times n}
    \end{array}\!\!\right]\nonumber = \bm A_0 + \mathbb{I}\Big(t\in [0,\tau]\Big)\,\bm{\delta  A}\,,
\end{align}
where to save the space we do not present explicit expressions for: the diagonal matrices ${\bm M}$ and ${\bm D}$ -- represent inertia and damping;  pre- and post-fault value ${\bm A}_0$ and its on-fault correction $\delta {\bm A}$ are due to the varying susceptance of the faulty line. Solution of (\ref{eq:swing-lin-matrix}) is explicitly given by
\begin{gather}\label{soldat}
\begin{aligned}
\bm x(t) &\!=\! e^{\tilde{\bm A} t}\bm x_0 \!+\! \int_0^t\!e^{\tilde {\bm A}(t-t_1)~} \!\bm P\,{\rm d}t_1\,, \forall t\le \tau \,,\\
\bm x(t) &\!=\! e^{\bm A_0 (t-\tau)}\bm x_\tau \!+\! \int_\tau^t\!e^{\bm A_0(t-t_1)~}\!\bm P\,{\rm d}t_1 \,, \forall t>\tau.\
\end{aligned}
\end{gather}
Where $\tilde{\bm A}=\bm A_0+\bm{\delta A}$. Eq.~\eqref{soldat} must usually be solved numerically. However, assuming that the initial condition $\bm x(t_0)$ and the eigendecompositions $\bm A = \bm U\, \bm \Lambda\, \bm U^{-1},\tilde{\bm A}=\tilde{\bm U}\tilde{\bm\Lambda} \tilde{\bm U}^{-1}$, with $\bm U\, \bm U^{-1}=\tilde{\bm U}\tilde{\bm U}^{-1}=\mathbb{I}$, are known, the system dynamics is governed by 
$$
\begin{aligned}
\dot{\bm \xi }& = \bm \Lambda\, \bm \xi + \bm U^{-1} \bm P\,, \;\forall t\le \tau \,,\nonumber\\ 
\dot{\bm \xi }& = \tilde{\bm \Lambda}\, \bm \xi + \tilde{\bm U}^{-1} \bm P\,,\;\forall t> \tau \,,\\ 
\bm \xi_0& = \bm U^{-1}\bm x_0\,,
\end{aligned}
$$
This initial value problem is piecewise constant for the faults we investigate in the present work.  It can be solved analytically for the eigenmodes $\bm \xi$, then converted into $\bm x$ as $\bm x(t) =\bm {\tilde U}\, \bm \xi(t),\ t\le \tau$ and $\bm x(t) =\bm U\, \bm \xi(t),\ t>\tau$, where tilde denotes that the is the eigenvectors of the perturbed matrix $\tilde{\bm A}$ that are used. More details can be found in \cite{pagnier2019optimal}.

\subsection{First Order Eigen Perturbation Expansion}

To simplify the computations, we assume that the perturbation induced by the fault is sufficiently small. This assumption is justified by the fact that the fault is both localized and short in duration ($\tau$). Under these conditions, a perturbation analysis (linear in the perturbation magnitude) is appropriate. We consider the perturbed matrix
$$
\tilde{\bm A}=\bm A_0+\bm{\delta A}=\bm A_0+ a\bm V,
$$
where $a$ is a dimensionless perturbation parameter and $\bm V$ is a perturbation matrix. Our goal is to approximate the spectrum of $\tilde{\bm A}$ in terms of the spectrum of the unperturbed matrix $\bm A_0$. If $a$ is small and the eigenvalues and eigenvectors of $\tilde{\bm A}$ are analytic functions of $a$ in a neighborhood of $\bar{0}$, the perturbed eigenvalue–eigenvector problem can be expressed as
\begin{gather}\label{evevprob}
\begin{aligned}
\tilde{\bm A}\bm {\tilde U}&=\bm {\tilde U}\tilde{\bm\Lambda},~\bm {\tilde U}^{-1}\tilde{\bm A}=\tilde{\bm\Lambda}\bm {\tilde U}^{-1},
\end{aligned}
\end{gather}
$$
\bm {\tilde U}=\sum_{i=0}^{\infty}a^i\bm {\tilde U}_i,~\bm {\tilde U}^{-1}=\sum_{i=0}^{\infty}a^i\left(\bm {\tilde U}^{-1}\right)_i,~\tilde{\bm\Lambda}=\sum_{i=0}^{\infty}a^i\tilde{\bm\Lambda}_i.
$$

Substituting these series into the eigenvalue equations (\ref{evevprob}) gives
$$
\begin{aligned}
\tilde{\bm A}\left[\sum_{i=0}^{\infty}a^i\bm {\tilde U}_i\right]&=\left[\sum_{i=0}^{\infty}a^i\bm {\tilde U}_i\right]\left[\sum_{i=0}^{\infty}a^i\tilde{\bm\Lambda}_i\right],\\
\left[\sum_{i=0}^{\infty}a^i\left(\bm {\tilde U}^{-1}\right)_i\right]&\tilde{\bm A}=\left[\sum_{i=0}^{\infty}a^i\tilde{\bm\Lambda}_i\right]\left[\sum_{i=0}^{\infty}a^i\left(\bm {\tilde U}^{-1}\right)_i\right].
\end{aligned}
$$
By identifying the zeroth-order terms with the unperturbed quantities, i.e., $\tilde{\bm\Lambda}_0 = \bm\Lambda$, $\bm {\tilde U}_0 = \bm U$, and $\left(\bm {\tilde U}^{-1}\right)_0 = \bm U^{-1}$, we can relate the perturbed spectrum of $\tilde{\bm A}$ to the spectrum of $\bm A_0$, the perturbation $\bm V$, and the parameter $a$. Following the methodology and using the same notation of \cite{bamieh2020tutorial}, the first-order approximation of the spectrum of $\tilde{\bm A} = \bm A_0 + a \bm V$ can be expressed as
\begin{gather}\label{ULUinv}
\begin{aligned}
\tilde{\bm A}&=\tilde{\bm U}\tilde{\bm\Lambda} \tilde{\bm U}^{-1}\approx \bar{\bm A}=\bar{\bm U}\bar{\bm\Lambda} \bar{\bm U}^{-1}\text{ : }\tilde{\bm A}=\bar{\bm A}+\mathcal{O}(a^2)\text{ and}\\
&\left[\bm \Pi^{+\circ}\right]_{ij}=\begin{cases}(\lambda_i-\lambda_j)^{-1}&:i\ne j\\0,&\text{otherwise}.\end{cases},\text{ yield the following}
\\
&\begin{cases}
\bar{\bm U}=\left(\bm U-a\bm U\left(\bm \Pi^{+\circ}\circ\left(\bm U^{-1}\bm V\bm U\right)\right) \right),\\
\bar{\bm\Lambda}=\left(\bm \Lambda+a\,\text{Diag}\left(\bm {U}^{-1}V\bm {U}\right)\right),\\
\bar{\bm U}^{-1}=\left(\bm U^{-1}+a\left(\bm \Pi^{+\circ}\circ\left(\bm U^{-1}\bm V\bm U\right)\right)\bm U^{-1}\right).
\end{cases}
\end{aligned}
\end{gather}
This way the \textbf{first order approximated solution} to (\ref{eq:swing-lin-matrix}) is explicitly given by
\begin{gather}\label{approx1}
\begin{aligned}
\bm x(t) &\!=\! \bar{\bm U}e^{\bar{\bm\Lambda}t}\bar{\bm U}^{-1}\bm x_0 \!+\! \int_0^t\!\bar{\bm U}e^{\bar{\bm\Lambda}(t-t_1)~} \!\bm\!\bar{\bm U}^{-1} \bm P\,{\rm d}t_1\,, \forall t\le \tau \,,\\
\bm x(t) &\!=\! e^{\bm A_0 (t-\tau)}\bm x_\tau \!+\! \int_\tau^t\!e^{\bm A_0(t-t_1)~}\!\bm P\,{\rm d}t_1 \,, \forall t>\tau.\
\end{aligned}
\end{gather}

This approach is both convenient and sufficiently accurate in the case of a \textbf{single-phase-to-ground fault} (see Section \textbf{\ref{Acc}}). However, it does not provide the same level of accuracy for a \textbf{three-phase fault}. In the following subsection, we present an efficient method to approximate solutions for the latter case.

\subsection{Eigenvalue Expansion in Multiple Steps}\label{Acc2}
Now, consider the case of a non-small localized perturbation with $a=1$, corresponding to a three-phase fault. As shown in Section \ref{Acc}, for small perturbations with $a:=m^{-1}\in(0,1)$, the error in the first-order approximation of the eigenpairs is bounded by $\mathcal{O}\big(a^2\big)$. This observation suggests that a perturbation of magnitude $a=1$ can be constructed by sequentially applying $m$ perturbations of size $a=\tfrac{1}{m}$.

The procedure is as follows. Starting from the unperturbed matrix $\bm A_0$, we approximate the spectrum of $\tilde{\bm A}=\bm A_0+\bm{\frac{1}{m}}V$ and compute $\bar{\bm A}=\bar{\bm U}\bar{\bm\Lambda} \bar{\bm U}^{-1}$. We then approximate the spectrum of $\bar{\bm A} + \frac{1}{m}\bm V$, and repeat this process iteratively $m$ times. After $m$ iterations, the $m$-th order approximation of the spectrum of $\tilde{\bm A} = \bm A_0 + \bm V$ is given by Eq. (\ref{ULUinv}) with $a = \frac{1}{m}$, where $(\bm U, \bm\Lambda)$ are taken from the previous iteration, i.e., from the first-order approximation of the spectrum at the $(m-2)$-th step. 

Upon completion of this iterative procedure, we obtain the triplet $\big(\bar{\bm U}^{\,(m)},\bar{\bm \Lambda}^{\,(m)}, \bar{\bm U}^{\,(m)-1}\big)$ which represents the $m$-th order approximation of the spectrum. The corresponding $m$-th order approximated solution to Eq. \eqref{eq:swing-lin-matrix} is therefore given by (\ref{approx1}) substituting $(\bar{\bm U},\bar{\bm \Lambda},\bar{\bm U}^{-1})$ for $\big(\bar   {\bm U}^{\,(m)},\bar{\bm \Lambda}^{\,(m)}, \bar{\bm U}^{\,(m)-1}\big)$.

\subsection{Accuracy}\label{Acc}

We now consider the perturbation corresponding to $a = \frac{1}{m}$, so that $\tilde{\bm A} = \bm A_0+\frac{1}{m}\bm V$. Since the perturbation is both small $\left(a = \tfrac{1}{m}\right)$ and localized ($\bm V$ rank-deficient), the result of \cite{eisenstat1998three} implies that the error between the first-order approximation $\bar{\lambda}_{\alpha}$ of the $\alpha$-th eigenvalue and its exact value $\tilde{\lambda}_{\alpha}$ satisfies the bound
\begin{equation}\nonumber\resizebox{.99999999999\hsize}{!}{$\big\|\tilde{\lambda}_{\alpha}-\bar{\lambda}_{\alpha}\big\|\le2a^2\left\|\bm V\right\|^2\big\|\big(\bm A-\lambda_{\alpha}\mathds{1}_{n\times n}\big)^{-1}\big\|\big\|[\bm U^{-1}]_{\alpha}\big\| \left\|[\bm U]^{\alpha}\right\|.$}
\end{equation}

Similarly, for the first-order approximations of the left and right eigenvectors, the following upper bound holds:
\begin{equation}\nonumber
{\begin{aligned}
\max&\Big\{\big\|\bm{\tilde{[U]}^{\alpha}}-\bm{\bar{[U]}^{\alpha}}\big\|,\big\|\bm{\tilde{[U^{-1}]}_{\alpha}}-\bm{\bar{[U^{-1}]}_{\alpha}}\big\|\Big\}\\
&\le \frac{2\left\|\bm V\right\|^2}{m^2}\big\|\big(\bm A-\lambda_{\alpha}\mathds{1}_{n\times n}\big)^{-1}\big\|^2\big\|[\bm U^{-1}]_{\alpha}\big\|\big\|[\bm U]^{\alpha}\big\|.
\end{aligned}}
\end{equation}
Using these $\mathcal{O}\big(m^{-2}\big)$, we can now estimate the error of the iterative scheme introduced in Section \textbf{\ref{Acc2}}:
\begin{equation}\nonumber
\small
{\begin{aligned}
\max&\Big\{\big\|\tilde{\bm\Lambda}-\bar{\bm \Lambda}^{(m)}\big\|,\big\|\tilde{\bm U}-\bar{\bm U}^{(m)}\big\|,\big\|\tilde{\bm U^{-1}}-\bar{\bm U}^{(m)-1}\big\|\Big\}\le\mathcal{O}\big(m^{-1}\big).
\end{aligned}}
\end{equation}
Here $\|\cdot\|$ refers to the spectral norm. Therefore, increasing the number of iterations (i.e., taking larger values of $m$) in the multi-step eigenvalue expansion improves the accuracy of the approximation, while the error decreases at the rate $\mathcal{O}\big(m^{-1}\big)$.

\subsection{Efficiency}

Analytical computation of first-order vector ODEs is often more efficient than numerical simulation when a closed-form solution exists. Such expressions provide exact results, eliminate iterative errors, and offer direct insight into system dynamics and parameter dependencies. In contrast, numerical methods—though more general—are computationally demanding and less interpretable, especially for stiff systems or long simulation horizons \cite{higham2008functions}. In this work, the matrices $\bm A$ and $\tilde{\bm A} = \bm A + a\bm V \in \mathbb{R}^{2n \times 2n}$ in Eq.~(\ref{soldat}) are assumed diagonalizable, allowing direct analytical computation. The procedure involves:
(i) eigendecomposition of $\tilde{\bm A} = \tilde{\bm U}\tilde{\bm \Lambda}\tilde{\bm U}^{-1}$ and $\bm A = \bm U\bm \Lambda\bm U^{-1}$, each costing $\mathcal{O}(n^3)$;
(ii) exponentiation of the diagonal matrices $\tilde{\bm \Lambda}$ and $\bm \Lambda$ at $\mathcal{O}(n)$; and
(iii) matrix–vector multiplications at $\mathcal{O}(n^2)$.
Thus, the total computational cost is $\mathcal{O}(n^3)$ for eigendecomposition and $\mathcal{O}(n^2)$ per evaluation of $\bm X(t)$. In contrast, numerical methods such as Runge–Kutta require $\mathcal{O}(n_{\Delta t}n^2)$ operations per evaluation. A performance comparison is given in Table~\ref{Theovssim}. The integral $S(\theta_{0\to T})$ in Eq.~(\ref{eq:over-heat-ij}) depends on both the system trajectory from Eq.~(\ref{soldat}) and the time-domain summation within the integrand. Assuming the interval $[0,T]$ is discretized into $n_{\Delta t} = T/\Delta t$ steps, each evaluation of $\bm X(t)$ costs $\mathcal{O}(n^2)$, leading to a total trajectory computation cost of $\mathcal{O}(n_{\Delta t}n^2)$. The summation over all transmission lines ${i,j} \in \mathcal{E}$ incurs $\mathcal{O}(|\mathcal{E}|)$ per time step, resulting in a total complexity of $\mathcal{O}(n_{\Delta t}|\mathcal{E}|)$. Hence, the overall computational cost for estimating $S(\theta_{0\to T})$ is $\mathcal{O}(n_{\Delta t}(n^2+|\mathcal{E}|))$. For sparse grids where $|\mathcal{E}| = \mathcal{O}(n)$, this simplifies to $\mathcal{O}(n_{\Delta t}n^2)$, while for dense grids where $|\mathcal{E}| = \mathcal{O}(n^2)$, the scaling remains $\mathcal{O}(n_{\Delta t}n^2)$ \cite{newman2018networks}. Spectral estimation via first-order perturbation requires one $\mathcal{O}(n^3)$ eigendecomposition of $\bm A$ and computation of the correction matrix $\Phi = \bm U^{-1}\bm V\bm U$ at $\mathcal{O}(n^2)$. Approximate eigenvalues for any perturbation parameter $a$ are then obtained from the diagonal of $\Phi$ at $\mathcal{O}(n)$. Therefore, both the exact spectrum of a single $\tilde{\bm A}$ and its first-order approximation have asymptotic cost $\mathcal{O}(n^2)$, but the approximation becomes significantly more efficient when multiple perturbations are considered, since preprocessing is performed once and each additional case requires only linear effort. Efficiency can be further improved through several strategies \cite{golub2013matrix}:
\begin{itemize}
\item \textbf{Sparsity:} exploit adjacency lists or sparse matrix formats to avoid redundant edge evaluations.
\item \textbf{Adaptive sampling:} concentrate computation on intervals where the indicator function is active.
\item \textbf{Edge pre-screening:} apply analytical bounds to exclude edges that never satisfy activation conditions.
\item \textbf{Parallelization:} distribute computations across time steps and edges using multi-core or GPU architectures.
\end{itemize}
\newpage

\section{Probabilistic Risk Assessment}
\label{sec:mcmc}

Assuming all the initial information from Section \textbf{\ref{sec:algorithm}}, we now focus on estimating the probability of exiting the safety polytope $\Pi_{\bm{\theta}}$ using the overload indicator $S(\theta_{0\to T})$ and an efficient sampling. Let $x_{0\to T}^{(0)}=x(0\to T)$ be the trajectory of the solution of the initial value problem stated in Eq.~\eqref{eq:swing-lin-matrix} and $S_{ij}\big(\theta^{(0)}_{0\to T}\big)$ be the line-specific overload indicator of $x^{(0)}_{0\to T}$ for all $\{i,j\}\in\cal E$:
\begin{equation}
\!\!S_{ij}\left(\theta^{(0)}_{0\to T}\right)
=\!\int\limits_0^T\!\!\mathbb{I} \bigg( \Big| \beta_{ij}(t)\left(\theta_i^{(0)}(t) \!-\! \theta_j^{(0)}(t)\right)\! \Big| \!>\! \bar{p}_{ij} \!\bigg){\rm d}t.
\end{equation}
We estimate the probability $Q_{ij}=\mathbb{P}\left[S_{ij}(\theta_{0\to T}) \ge \gamma\right]=\mathbb{E}\left[\mathbb{I}_{\{S_{ij}(\theta_{0\to T}) \ge \gamma\}}\right]$ for a fixed threshold $\gamma > 0$, i.e. \textbf{the likelihood of overloading line $\{i,j\}\in\cal E$ for more than $\gamma$ seconds}. An efficient and accurate way to approximate $Q_{ij}$ is by using \textbf{adaptive importance sampling} via the \textbf{cross entropy method} \cite{rubinstein_optimization_1997,de2005tutorial}, we briefly explain how does this work and how do we implement it for this particular case:

\subsection{Cross-Entropy Method (CEM)}\label{subsec:cem}
CEM \cite{rubinstein_optimization_1997,de2005tutorial} adaptively fits an importance law to concentrate samples in the rare-event region
\[
\mathcal{R}_\gamma \;=\;\bigl\{Z=(\alpha,\tau)\in \left(\{\varnothing\}\cup\mathcal E\right) \;\times\; \mathbb{R}_+ \;\big|\; S_{ij}(\theta_{0\to T})\ge \gamma\bigr\}.
\]
We model via a simple parametric family characterized by the sampling density:
$q_\nu(\alpha,\tau)=\phi_\alpha\,\mathrm{Pois}(\tau;\lambda_\alpha)$, where $\phi_\alpha\geq 0,\ \alpha\in{\cal E}$, $\sum_{\alpha\in {\cal E}}\phi_\alpha=1$, $\tau\in \mathbb{R}_+$ and with the parameter $\nu\doteq (\phi_\alpha,\lambda_\alpha|\alpha\in{\cal E})$.
Starting from $\nu^{(0)}$, at iteration $t$ we draw $Z_k\sim q_{\nu^{(t-1)}}$ samples, compute trajectories $\theta^{(k)}_{0\to T}$ and importance weights

$$
w_k \propto \mathbf 1\{S_{ij}(\theta^{(k)}_{0\to T})\ge\gamma\}\,\frac{p_Z(Z_k)}{q_{\nu^{(t-1)}}(Z_k)}
$$ 

also restricting to an ``elite'' subset. The CE update maximizes the weighted log–likelihood. Iterating until stabilization gives $q_{Z}^\star=q_{\nu^\star}$. The IS estimator of \(\;Q_{ij}=\mathbb P_{p_Z}\big[S_{ij}(\theta_{0\to T})\ge\gamma\big]\;\) is
\begin{gather*}
\hat{Q}_{ij}
=\frac{1}{N}\sum_{k=1}^N
\mathbf 1\!\left(S_{ij}(\theta^{(k)}_{0\to T})\ge\gamma\right)
\frac{p_Z(\alpha_k,\tau_k)}{q_{Z}^\star(\alpha_k,\tau_k)}\,,
\ (\alpha_k,\tau_k)\sim q_{Z}^\star.
\end{gather*}

\subsection{N1Plus}\label{subsec:n1plus}
We introduce N1Plus, a computational engine designed to efficiently simulate the dynamic response of transmission power grids to fault events, capturing system behavior across sub-second to minute timescales. The model perturbs parameters in a linearized form of the swing equations, keeping them fixed prior to transient stabilization. Fault events are sampled, and the resulting system dynamics are resolved analytically using the closed-form solutions in Eqs.~\eqref{soldat}, avoiding full nonlinear simulations. Based on these simulations, N1Plus estimates the probability of exiting the safety polytope $\Pi_{\bm{\theta}}$ through the cumulative overload indicator $S(\theta_{0\to T})$, as well as line-specific indicators that identify critical network components. To evaluate the probability of extreme events, stochastic sampling techniques are used to approximate the distribution of line overloads triggered by random disturbances. This framework represents a dynamic and stochastic generalization of the classical $N-1$ contingency analysis, providing a novel tool for operational resilience assessment. Conceptually, it extends the static instanton methodology proposed in \cite{chertkov2010predicting} and builds upon the rare-event sampling techniques developed in \cite{owen2019importance}. The methodology underlying N1Plus is summarized in Algorithm~\ref{alg:n1plus}.

\begin{algorithm}
  \caption{N1Plus Dynamic Kernel (pseudo‑code)}
  \label{alg:n1plus}
  \begin{algorithmic}[1] 
  \Require ${\cal G}=({\cal V},{\cal E}),L,M,D,P,\Pi_{\bm{\theta}}$ (\text{Section }\textbf{\ref{sec:algorithm}}), 
  $\lambda\in [0.1,0.9]$ \textit{rate of fault's duration}, $T>0$ \textit{duration of simulation}, 
  $N>0$ \textit{sample size}, $\gamma\ge 0$ \textit{threshold}, 
  parametric family $q_Z(\cdot;\nu)$ with initial parameters $\nu^{(0)}$.
  \State Initialize iteration counter $t\gets 0$.
  \Repeat
      \State $t\gets t+1$
      \For{$i=1,\dots,N$}
          \State $(\alpha_i,\tau_i)\sim q_Z\left(\cdot;\nu^{(t-1)}\right)$ \Comment{Fault dur.+tripped line.}
          \State Solve~(\ref{eq:swing-lin-matrix}) for $(\alpha_i,\tau_i)$: $\theta^{(i)}_{0\to T}$ \Comment{Use \ref{Acc2}.}
          \State $S_i\gets S\!\left(\theta^{(i)}_{0\to T}\right)$ \Comment{Use (\ref{glonOHI}).}
      \EndFor
      \State $E^{(t)} = \{ (\alpha_i,\tau_i): S_i \ge \gamma^{(t)}\}$\Comment{Elite set.}
      \State $\gamma^{(t)}:(1-\rho)$-quantile of $\{S_i\}$\Comment{$\rho\in(0,1)$.}
      \State $\nu^{(t)} = 
      \underset{\nu}{\mathrm{argmax}}\left\{\sum_{(\alpha_i,\tau_i)\in E^{(t)}}\
      \log \left(q_Z(\alpha_i,\tau_i;\nu)\right)\right\}$
  \Until{Convergence of $\nu^{(t)}$: optimal distribution $q_Z^\star\left(\cdot;\nu^{(t)}\right)$.}
  \State Draw $\{(\alpha_k,\tau_k)\sim q_Z^\star\}_{k=1}^N$ from $q_Z^\star\left(\cdot;\nu^{(t)}\right)$.
  \State For each $(\alpha_k,\tau_k)$, solve~(\ref{eq:swing-lin-matrix}) and compute $S_{lm}\!\left(\theta^{(k)}_{0\to T}\right)$.
  \State Estimate $\hat{Q}_{lm}$ \Comment{Use (\ref{subsec:cem}).}
  \Return $\hat{Q}\doteq\left(\hat{Q}_{lm}\right)_{\{l,m\}\in{\cal E}}$ \Comment{Repeat 14-15 $\forall\{l,m\}\in{\cal E}$.}
  \end{algorithmic}
\end{algorithm}

\section{Case Study: Israeli Transmission Grid}
\label{sec:casestudy}

We apply the \textbf{N1Plus Algorithm~\ref{alg:n1plus}} to an open-source model of the Israeli Grid (IG), consisting of 32 buses (11 generators and 21 loads) and 36 transmission lines. Details of this model, along with the parameters required by Algorithm~\ref{alg:n1plus} and Section~\textbf{\ref{sec:algorithm}}, can be found in \cite{israelinfo}. Our evaluation focuses on both the accuracy and computational efficiency of the proposed framework. The final implementation will be released as a Julia-based open-source toolkit, tentatively named \textbf{N1Plus.jl}. Sections~\textbf{\ref{sec:3PF}} and~\textbf{\ref{sec:1PF}} present numerical results for single-phase and three-phase fault scenarios, respectively. In these studies, we analyze the statistical behavior of the line-specific overload indicator $S_{ij}(\theta_{0\to T})$ as a function of the random fault duration $\tau$, across a range of single- and three-phase disturbances. Specifically, we examine the probability $Q_{ij}(\tau)=\mathbb{P}\left[S_{ij}\geq T^{(*)}\right]$, which measures the likelihood that the overload indicator exceeds a predefined threshold $T^{(*)}$ for a given fault duration. Thresholds of $T^{(*)}=5\text{s}$ for single-phase faults and $T^{(*)}=0.125\text{s}$ for three-phase faults are adopted. To interpret event severity, we employ a qualitative risk classification based on $Q_{ij}$: values above $10\%$ are labeled red (critical risk), those between $5\%$ and $10\%$ as yellow (moderate risk), and those below $5\%$ as green (low risk). These thresholds are illustrative and not intended to reflect operational standards. In practical settings, such values should be defined by system operators using empirical data, engineering judgment, and domain-specific expertise. Finally, Section~\textbf{\ref{BenchAn}} presents the comparison, validation, and performance analysis of both the explicit and approximated solutions of Eq.~\eqref{eq:swing-lin-matrix}. This section highlights the conditions under which it is advantageous—or not—to use the approximation in Eq.~\eqref{Acc2} instead of the exact formulation in Eq.~\eqref{soldat}. A comprehensive relative error analysis is conducted over all possible random line trippings with varying magnitudes, together with a detailed evaluation of computation times.

\subsection{Statistics of the Three Phase Fault with $T^{(*)}=5\text{s}$}\label{sec:3PF}

\begin{figure}[thb]
    \centering
	\includegraphics[width=1.0\columnwidth]{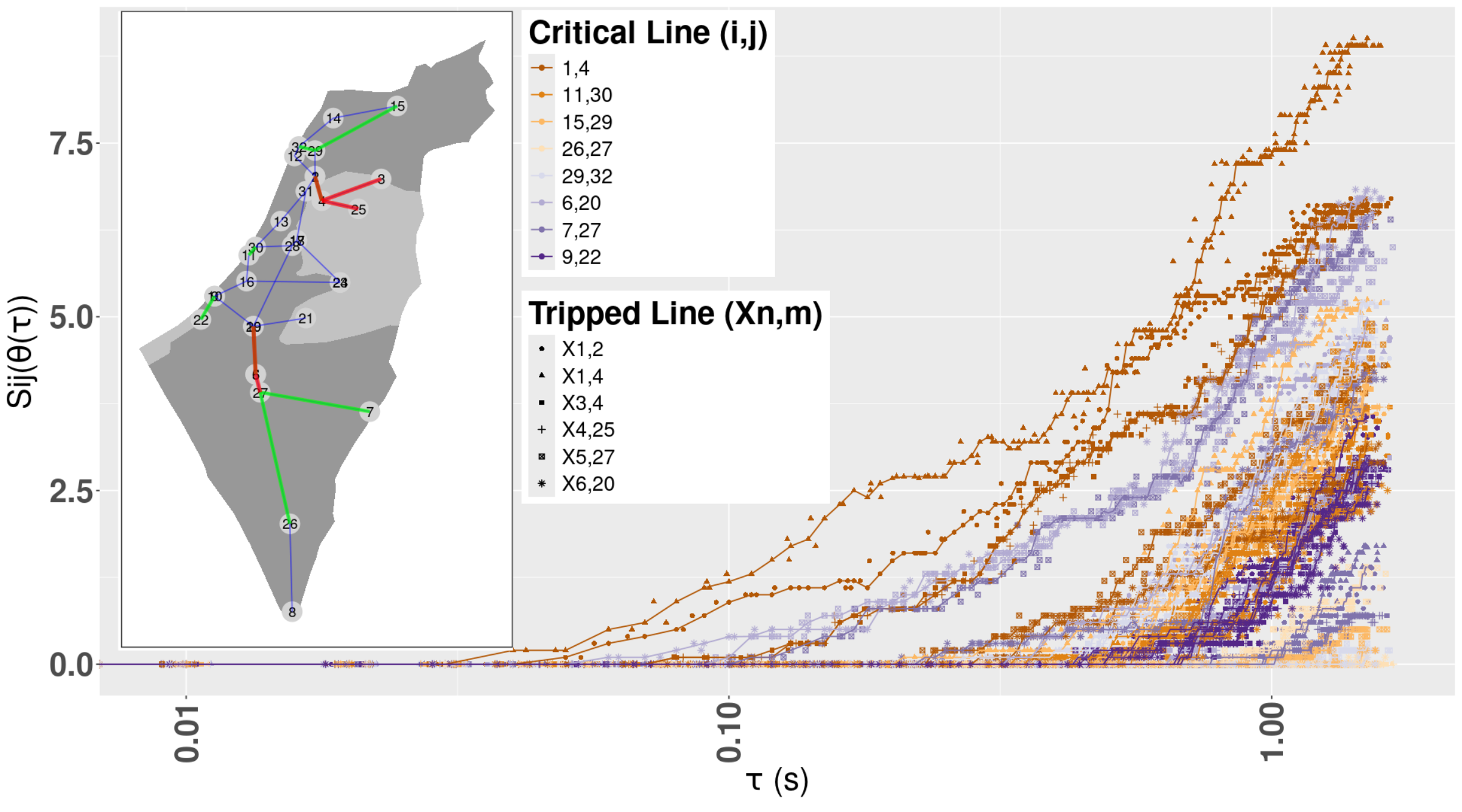}
	\caption{Overload indicator vs. fault duration (3-Phase Faults)· Each curve shows line-specific overload vs. fault time $\tau$. \textbf{\textcolor{red}{Tripped lines}} $(X_{n,m})$: trigger peak overload \textbf{\textcolor{OliveGreen}{Critical lines} (associated with transformers)} $(i,j)$: \textbf{accumulate} highest total overload. Left: Israel grid with \textbf{relevant lines} highlighted.}
	\label{fig2}
\end{figure}

\begin{figure}[thb]
    \centering
	\includegraphics[width=1.0\columnwidth]{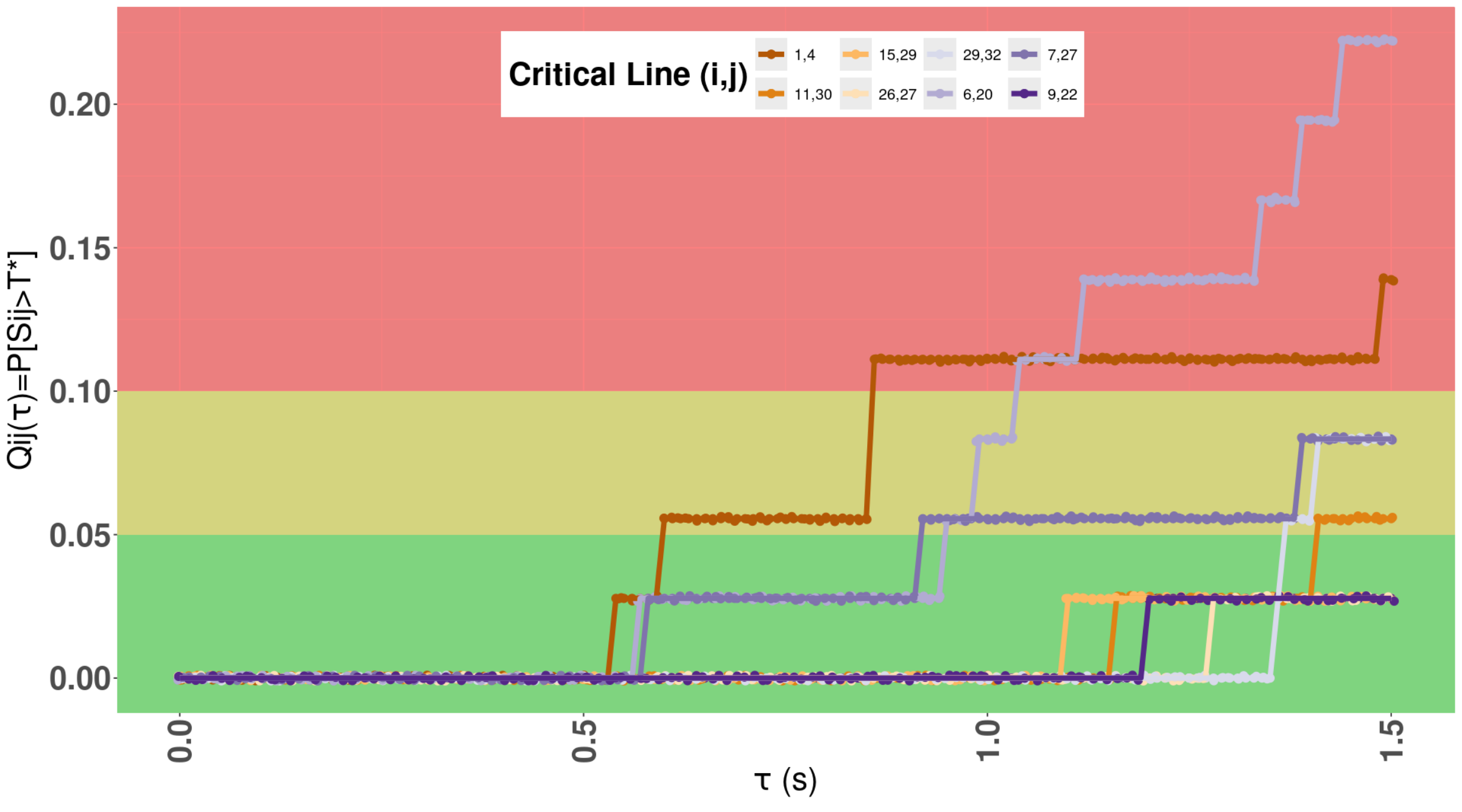}
	\caption{Probabilistic security assessment (3-Phase Faults). For critical lines and fault durations $\tau$ we determine $\mathbb{P}\left[S_{ij}\ge T^{(*)}\right]$. Risk zones by $Q_{ij}(\tau)<5\%$: \textbf{\textcolor{OliveGreen}{Safety Zone}}, $Q_{ij}(\tau)\in [5\%,10\%]$: \textbf{\textcolor{YellowOrange}{Warning Zone}} and $Q_{ij}(\tau)>10\%$: \textbf{\textcolor{LightRed}{Emergency Zone}}}
	\label{fig3}
\end{figure}

Figs.~(\ref{fig2}) and (\ref{fig3}) present the results obtained from a total of $N = 100{,}000$ simulations of system of Eqs.~\eqref{eq:swing-lin-matrix}. In these simulations, each fault corresponds to the complete removal of a transmission line. The tripped line in each sample is selected uniformly at random from all lines in the grid, and the fault duration is drawn from an exponential distribution with rate parameter $\lambda = 0.1$.

\subsection{Statistics of the Single Phase Fault with $T^{(*)}=0.125\text{s}$}\label{sec:1PF}

\begin{figure}[thb]
    \centering
	\includegraphics[width=1.0\columnwidth]{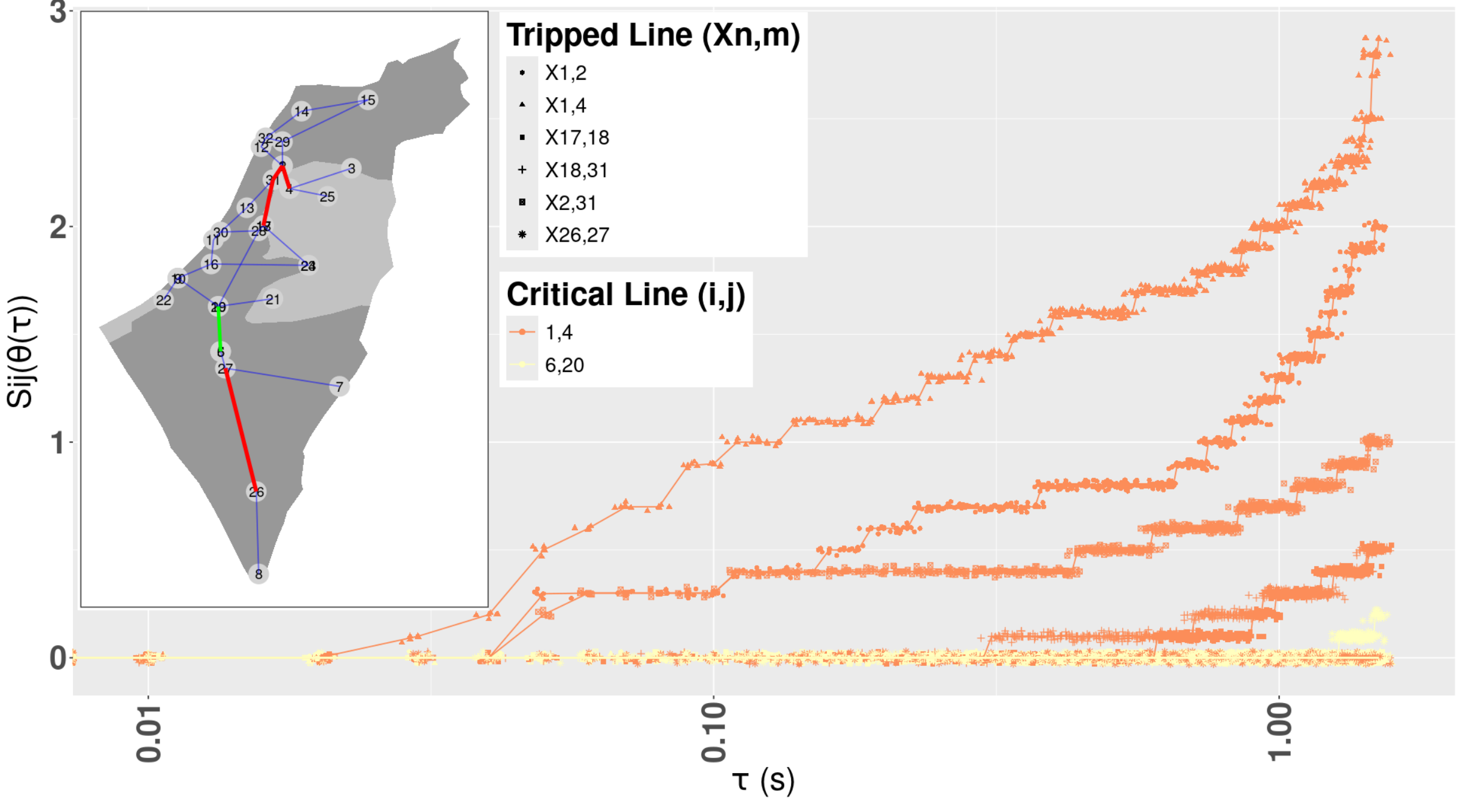}
	\caption{Overload indicator vs. fault duration for the case of single phase faults, otherwise notations and descriptions are the same as in Fig.~(\ref{fig2}).}
	\label{fig4}       
\end{figure}

\begin{figure}[thb]
    \centering
	\includegraphics[width=1.0\columnwidth]{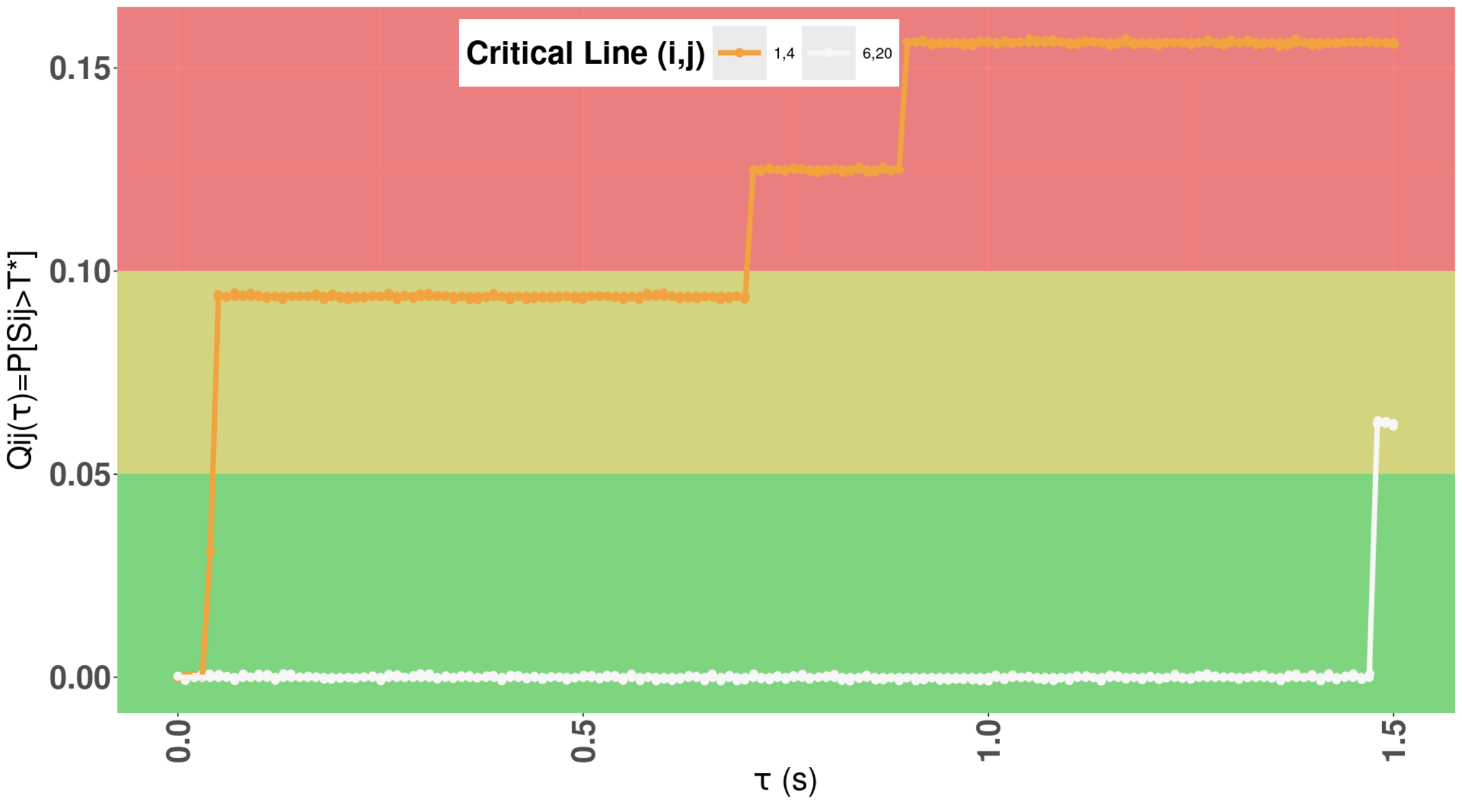}
	\caption{Probabilistic security assessment for the case of single phase faults and critical lines shown in Fig. (\ref{fig4}), otherwise notations and descriptions are the same as in Fig.~(\ref{fig3}).}
    \label{fig5}
\end{figure}

Figs.~(\ref{fig4}) and (\ref{fig5}) present the results obtained from a total of $N = 100{,}000$ simulations of the system of Eqs.~\eqref{eq:swing-lin-matrix}. In this case, each fault corresponds to a partial removal of a transmission line, modeled as a one-third reduction in line susceptance (see Footnote \ref{footnote:disclaimer}). The tripped line in each simulation is selected uniformly at random from all lines in the grid, and the fault duration is drawn from an exponential distribution with rate parameter $\lambda = 0.1$.

\subsection{Benchmark Analysis}\label{BenchAn}

\begin{table*}[!t]
\caption{Comparison in Processing Time: Exact Solution (\textbf{\ref{soldat}}) vs Approx Solution (\textbf{\ref{Acc2}}) in \textbf{m} steps. \textbf{{\color[HTML]{3531FF} \textbf{Blue}}}: Avg. computation time + std. dev $(\sigma)$ of (\ref{soldat}).
      \textbf{{\color[HTML]{009901} \textbf{Green}}}: Approximation (Sec. \ref{Acc2}) faster than exact.
      \textbf{{\color[HTML]{CB0000} \textbf{Red}}}: Slower than exact.}
  \label{proctime}
\small
\centering
\begin{tabular}{lrr}
\specialrule{.1em}{.05em}{.05em}
\textbf{Type of Solution} & \multicolumn{1}{c}{\textbf{Mean Processing Time [s]}} & \multicolumn{1}{c}{\textbf{Standard Deviation}} \vphantom{$\int_p^f$}                          \\ \specialrule{.1em}{.05em}{.05em}
{Exact Solution}   & {\color[HTML]{3531FF} \textbf{0.815}} & {\color[HTML]{3531FF} \textbf{0.186}}  \\ 
{m = 10 (Approx)}    & {\color[HTML]{009901} \textbf{0.025}} & {\color[HTML]{009901} \textbf{0.027}} \\
{m = 40 (Approx)}    & {\color[HTML]{009901} \textbf{0.458}} & {\color[HTML]{009901} \textbf{0.055}} \\ 
{m = 70 (Approx)}    & {\color[HTML]{009901} \textbf{0.661}} & {\color[HTML]{009901} \textbf{0.063}} \\ 
{m = 100 (Approx)}   & {\color[HTML]{009901} \textbf{0.811}} & {\color[HTML]{009901} \textbf{0.077}} \\ 
{m = 150 (Approx)}   & {\color[HTML]{CB0000} \textbf{1.129}} & {\color[HTML]{CB0000} \textbf{0.111}} \\ 
{m = 200 (Approx)}   & {\color[HTML]{CB0000} \textbf{1.520}} & {\color[HTML]{CB0000} \textbf{0.123}} \\ \hline
\end{tabular}
\end{table*}

\begin{table*}[!t]
\caption{Processing Time \& Convergence Comparison: Cross-Entropy (\ref{subsec:cem}) vs. Monte Carlo \cite{dagum2000optimal} for global overload (\ref{glonOHI}) probability at thresholds \textbf{$(\gamma)$}. \textbf{\textcolor[HTML]{009901}{CEM estimate}} converges faster (time \& sample size) than \textbf{\textcolor[HTML]{CB0000}{MC estimate}}.}
  \label{table4}
\small
\centering
\begin{tabular}{rrrrrrr}
\hline
 & \multicolumn{3}{c}{\textbf{\textcolor[HTML]{009901}{Cross Entropy Method} }} &
\multicolumn{3}{c}{\textbf{\textcolor[HTML]{CB0000}{Monte Carlo Method}}}
\\
\multicolumn{1}{c}{$\gamma$ [s]} &\multicolumn{1}{c}{Sample Size} & \multicolumn{1}{c}{$\hat{Q}\approx \mathbb{P}[S > \gamma]$}  & \multicolumn{1}{c}{Mean Time [s]} & \multicolumn{1}{c}{Sample Size} & \multicolumn{1}{c}{$\hat{Q}\approx \mathbb{P}[S > \gamma]$} &\multicolumn{1}{c}{ Mean Time [s]}\vphantom{$\int_p^f$}\\
\hline
\vphantom{$\int_p^f$}0.5 & 2500  & {\color[HTML]{009901}\textbf{0.07538}} & {\color[HTML]{009901}\textbf{\textit{$1.656\times 10^3$}}} &  25000 & {\color[HTML]{CB0000}\textbf{0.07401}} & {\color[HTML]{CB0000}\textbf{\textit{$1.792\times 10^4$}}} \\
5.0 & 5000  & {\color[HTML]{009901}\textbf{0.00319}} & {\color[HTML]{009901}\textbf{\textit{$3.312\times 10^3$}}} & 25000 & {\color[HTML]{CB0000}\textbf{0.00319}} & {\color[HTML]{CB0000}\textbf{\textit{$1.792\times 10^4$}}} \\
10.0 & 10000  & {\color[HTML]{009901}\textbf{0.00042}} & {\color[HTML]{009901}\textbf{\textit{$6.666\times 10^3$}}} & 50000 & {\color[HTML]{CB0000}\textbf{0.00041}} & {\color[HTML]{CB0000}\textbf{\textit{$3.585\times 10^4$}}} \\
\hline
\end{tabular}
\end{table*}

\begin{figure}[thb]
    \centering
	\includegraphics[width=1.0\columnwidth]{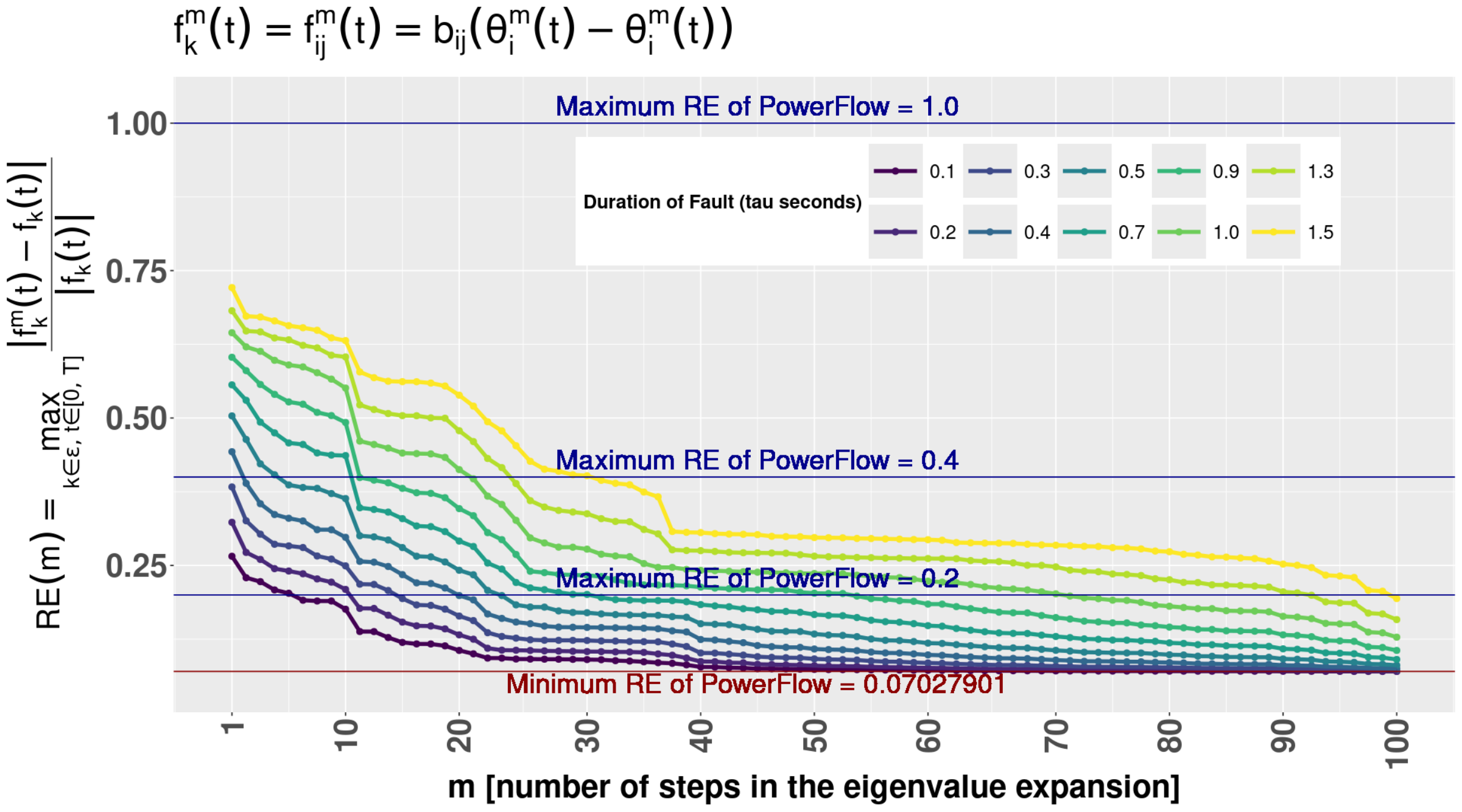}
	\caption{Maximum relative error between the exact (\ref{soldat}) and approximated solutions for all tripped lines, with fault duration in the interval $\tau \in \left[0s,1.5 \text{s}\right]$ and approximation steps (\ref{Acc2}) parameterized by $m \in [0,100]$.}
	\label{fig6}
\end{figure}

Figure \ref{fig6} and Table \ref{proctime} summarize the benchmark results obtained from solving system (\ref{eq:swing-lin-matrix}) over a total of $N=100{,}000$ simulations. First, we evaluate the optimal number of steps $m$ in the iterative scheme (\ref{Acc2}) that minimizes the maximum recorded relative error across all possible faults in the IG power grid for different fault durations $\tau$. We then compare the average processing times between the exact and approximated solutions. Consistent with Section \ref{sec:3PF}, we restrict the analysis to scenarios involving complete line trippings, as these faults induce the largest oscillations in the system dynamics over time.

Finally, Table~\ref{table4} summarizes the efficiency of estimating the probability of global system overloading using the CEM approach (Section~\ref{subsec:cem}). Analogous to Section~\ref{sec:mcmc}, we approximate $Q=\mathbb{P}\left[S(\theta_{0\to T}) \ge \gamma\right]$ for threshold levels $\gamma \in \{0.5, 5, 10\}$ and sample sizes $N \in \{1000, 2500, 5000, 10000, 25000, 50000\}$. The table highlights the minimum sample size required for convergence and reports the corresponding mean processing time per run. These results are compared against those obtained using traditional Monte Carlo sampling \cite{dagum2000optimal}, demonstrating the superior convergence efficiency of the CEM estimator.

\section{Conclusion and Outlook} 
\label{sec:conclusion}
We have developed a \textit{dynamic yet analytically tractable} framework for real-time \mbox{$N-1$} security assessment. Starting from a linearized swing-equation formulation, we derived closed-form solutions that bypass costly time-domain simulations, introduced an \textbf{Overload Indicator} to quantify transient severity in a risk-consistent manner, and implemented an efficient sampling engine -- \textbf{N1\textsc{Plus}} -- that transforms these analytics into an operator-ready advisory tool. Validation on the Israeli transmission grid demonstrated (i) sub-second evaluation times on a standard laptop for $\mathcal{O}(10^5)$ fault trajectories, and (ii) excellent agreement with benchmark simulations, even for severe single-phase faults.

\noindent\textbf{Next steps.}
Future development of \textsc{N1Plus} will focus on:
\begin{enumerate}
\item \textbf{Scalability and sparsity.} Porting analytics to sparse GPU kernels and hierarchical solvers to enable screening of multi-area systems ($n\sim 10^4$ buses) within the SCADA/EMS cycle.

\item \textbf{Nonlinear and voltage-coupled dynamics.} Incorporating reduced-order voltage and converter models to capture minute-scale voltage excursions and low-inertia effects (missing in the current version).

\item \textbf{Data-driven parameter tuning.} Integrating PMU/SCADA data with physics-informed priors via Bayesian updating for online calibration of fault statistics and protection clearing times.

\item \textbf{Uncertainty quantification.} Formalizing safety polytope “recurrence” through return-time distributions under stochastic injections and component uncertainty.

\end{enumerate}

\noindent\textbf{Operational and Planning Value.} Because \textsc{N1Plus} maintains analytical transparency, it can serve as both a real-time monitoring dashboard and a strategic planning tool. In day-ahead or seasonal studies, it can rank latent vulnerabilities, and propose targeted reinforcements by quantifying the probabilities of rare contingencies. In light of the April 2025 Iberian blackout \cite{entsoe2025iberian} -- where cascading failures and voltage instability stressed conventional static screening methods \textsc{N1Plus} offers a complementary dynamic perspective. By estimating the risk and trajectories of transient overloads, it helps operators detect hidden failure paths that static $N-1$ analyses may miss. Ultimately, \textsc{N1Plus} bridges the divide between high–fidelity EMT simulators (accurate but slow) and static $N-1$ screens (fast but dynamics-blind), furnishing operators with a principled, efficient tool for assessing transient risk in evolving power systems -- especially in grids under high renewable penetration and low inertia such as Spain’s.

\end{document}